\documentclass{elsart}

\usepackage{natbib}

\usepackage{amsmath}
\usepackage{amsfonts}

\newtheorem{theorem}{Theorem}[section]
\newtheorem{lemma}[theorem]{Lemma}
\newtheorem{proposition}[theorem]{Proposition}
\newtheorem{corollary}[theorem]{Corollary}
\newtheorem{example}[theorem]{Example}
\newtheorem{remark}[theorem]{Remark}
\newenvironment{pot}{\noindent{\it Proof of Theorem \ref{main}.}}{\hfill$\Box$}
\newenvironment{proof}{\noindent{\bf Proof.}}{\hfill$\Box$}

\def\AA{{\mathcal A}}
\def\EE{{\mathcal E}}
\def\rr{{\mathcal R}}
\def\BB{{\mathcal B}}
\def\MM{{\mathbb M}}
\def\KK{{\mathbb K}}
\def\PP{{\mathbb P}}
\def\NN{{\mathbb N}}
\def\QQ{{\mathbb Q}}
\def\ZZ{{\mathbb Z}}

\begin{document}

\begin{frontmatter}

\title{Subresultants and Generic Monomial Bases}
\author[D'Andrea]{Carlos D'Andrea}
\author[Jeronimo]{Gabriela Jeronimo}

\address[D'Andrea]{Miller Institute for Basic Research in Science and Department of Mathematics
University of California, Berkeley CA 94720-3840, USA.
{\tt cdandrea@math.berkeley.edu}}

\address[Jeronimo]{Departamento de Matem\'atica, Facultad de Ciencias Exactas y Naturales, Universidad de Buenos Aires.\\
Ciudad Universitaria, Pabell\'on I, 1428 Buenos Aires, Argentina.\\
{\tt jeronimo@dm.uba.ar}}


\begin{abstract}
Given $n$ polynomials in $n$ variables of respective degrees
$d_1,\dots,d_n,$
and a set of monomials of cardinality $d_1\dots d_n,$ we give an explicit
subresultant-based polynomial expression in the coefficients of the input polynomials whose non-vanishing is a
necessary and sufficient condition for this set of monomials to be a basis
of the ring of polynomials in $n$ variables modulo the ideal generated by the system of polynomials. This approach allows
us to clarify the algorithms for the B\'ezout
construction of the resultant.
\end{abstract}
\end{frontmatter}

\section{Introduction}\label{uno}
Consider a system of $n$ polynomials in $n$ variables with
coefficients in a field $\KK,$
$f_1(x_1,\dots,x_n),\dots,f_n(x_1,\dots,x_n),$ with respective
degrees $d_1,\dots,d_n$. Generically, this system has ${\mathbf
d}:= d_1.d_2.\dots .d_n$ roots in the algebraic closure of $\KK.$
This is the very well-known {\it B\'ezout formula} which appeared
in \citet{Bez} (see \citet{CLO1} for a modern treatment of this).
\par
One can say something more about what ``generic'' means above: let $V(f_1,\dots,f_n)\subset
\overline{\KK}^n$ be the set of common zeros of the polynomials $f_1,\dots,f_n,$ and set
$$f_i:=\sum_{j=0}^{d_i} f_{ij},\qquad i=1,\dots,n,
$$
where $f_{ij}$ is the homogeneous component of $f_i$ of degree $j.$ Then,
it turns out that $V(f_1,\dots,f_n)$ is a finite set and its cardinality (counting multiplicities)
is ${\bf d}$ if and only if the system of homogeneous equations
\begin{equation}\label{infty}
f_{1d_1}=0,\,f_{2d_2}=0,\,\dots,f_{nd_n}=0
\end{equation}
has no solution in projective space $\PP^{n-1}$---see \citep[Ch.
3, Thm. 5.5]{CLO2} for a proof of this result and also \citep[Ch.
4, Definition 2.1]{CLO2} for the definition of multiplicity of a
zero of a polynomial system.
\par
{}From a more algebraic point of view, if we set $I:=(
f_1,\dots,f_n)$ for the ideal generated by the $f_i$'s in
$\KK[x_1,\dots,x_n],$ the fact that $V(I)\subset\overline{\KK}^n$
has ${\bf d}$ points counted with multiplicity  means that the
${\KK}$-algebra $\AA:={\KK}[x_1,\dots,x_n] / I$ is a
${\KK}$-vector space of dimension  ${\bf d}$. As $\AA$ is
generated by the set of (the images in $\AA$ of) all monomials in
${\KK}[x_1,\dots,x_n],$ one can always find a basis of monomials
for $\AA$ (finite or not).
\par
In this paper, we will focus our attention on the following problem: given a set $\MM$ of
${\bf d}$ monomials, how can we decide if they are a basis of $\AA$ or not?
\par
We could use Gr\"obner bases for solving this problem, but we
would like our answer to be a function on the input set $\MM$
only, and not depending on an extra monomial ordering and other
intermediate steps that are needed in Gr\"obner bases algorithms.
\par
One of the main results of this paper is a polynomial expression
in the coefficients of $f_1,\dots,f_n$ which vanishes if and only
if the set $\MM$ fails to be a basis of $\AA.$ The expression we
get can be described in terms of resultants and subresultants of
homogeneous polynomials obtained from the input system, which is
the algebraic counterpart of this problem in the homogeneous case
 \citep[see][{}]{CLO2,Cha1,sza}.
\par
The problem of deciding whether a given set of monomials $\MM$ is a basis of $\AA$ or not is important in elimination theory due to the fact that algorithms for computing
resultants, B\'ezout identities, reduction modulo an ideal and explicit versions of the Shape Lemma can be reduced to linear algebra computations in the quotient ring, avoiding the use of
Gr\"obner bases, if one succeeds in finding such a basis $\MM$.
\par
\citet{Bez} was the first to work following this approach, which
was extended by  \citet{Mac}, who answered this question in the
case $ \MM=\{x_1^{\alpha_1}\dots x_n^{\alpha_n},\,0\leq
\alpha_i\leq d_i-1\} $ by means of a polynomial expression in the
coefficients of the input polynomials \citep[see also][{}]{Mac2}.
Our results, when applied to Macaulay's case, recover his
original formulation.
\par
In this direction, some results were obtained by \citet{Cha3},
provided that all the $f_i$'s are generic and homogeneous. If the
input system is generic and sparse, a generalization of the case
we are dealing with here, partial results were obtained by
\citet{ER} and \citet{PS} for $\MM$'s constructed by means of
regular triangulations of polytopes.
\par
A different approach based on recursive linear algebra is
provided in \citet{BU} for specific $\MM.$ In Section \ref{seis},
we will compare our results with those obtained in this article.
\par
The paper is organized as follows: some preliminary results are
stated in Section \ref{dos}. In Section \ref{pretres}, we recall
the definition and basic properties of multivariate
subresultants, as introduced in \citet{Cha1}. We relate
subresultants with our problem in Section \ref{tres}, associating
with any given set $\MM$ a polynomial whose non vanishing is
equivalent to the fact that $\MM$ is a basis of $\AA$. In Section
\ref{cuatro}, we show that, for certain $\MM$'s, this polynomial
expression depends only on the coefficients of
$f_{1d_1},\dots,f_{nd_n},$ and moreover, it can be decomposed into
factors. Then, we give in Section \ref{cinco} some rational
expressions for generalized Vandermonde determinants. These
results, along with those presented in Section \ref{cuatro},
allow us a better understanding of the recursive algorithm
proposed in \citet{BU}. Finally, we conclude by comparing our
results with those obtained in \citet{BU} in Section \ref{seis}.

\section{Preliminary Results}\label{dos}

Let ${\rm Res}_{d_1,\dots,d_n} (^{\bf .})$ be the homogeneous
resultant operator, as defined in \citet{Mac,VdW,CLO2}. We recall
the following well-known result \citep[see][for a proof]{CLO2}:

\begin{proposition}\label{affine}
The system (\ref{infty}) has a nontrivial solution in $\overline{\KK}^n$ if and only if
${\rm Res}_{d_1,\dots,d_n}(f_{1d_1},\dots,f_{nd_n})=0.$
\end{proposition}

\medskip
\begin{remark}
This proposition, together with our previous remarks about the
quotient ring $\AA$, gives a proof for the Choice Conjecture
stated in \citet{BU}: The condition ${\rm
Res}_{d_1,\dots,d_n}(f_{1d_1},\dots,f_{nd_n})\ne 0$ is necessary
and sufficient for the existence of a set $\MM$ of ${\bf d}$
monomials which is a basis of $\AA$ (and hence, any polynomial
can be reduced with respect to this set). Of course, the hard
problem is to find such an $\MM$!
\end{remark}

\par
Let $\KK$ be a field, $f_1,\dots,f_n\in\KK[x_1,\dots,x_n]$ and
$$\MM:=\{m_1,\dots,m_{\bf d}\}\subset\KK[x_1,\dots,x_n]$$ be a set of ${\bf d}$ monomials.
Set $\rho := d_1+\dots+d_n-n,$ and $$\delta :=
\delta(\MM)=\max\{\deg(m_i),\,i=1,\dots,{\bf d}\}.
$$
Let $x_0$ be a new variable. For every polynomial
$p(x_1,\dots,x_n)\in\KK[x_1,\dots,x_n]$ we define
$$p^0(x_0,x_1,\dots,x_n):=x_0^{\deg(p)}p\big(\frac{x_1}{x_0},\dots,\frac{x_n}{x_0}\big),$$
i.e. $p^0$ is the homogenization of $p$ with a new variable
$x_0$, and for every $t\geq\delta,$ we set
$$\MM_t:=\{mx_0^{t-\deg(m)},\,m\in\MM\}.$$

Let $\AA_0$ be the quotient ring ${\KK}[x_0,\dots,x_n] /
(f^0_1,\dots,f^0_n).$ It is a graded ring of the form
$\AA_0=\bigoplus_{i=0}^\infty \AA_{0i}.$

Set $H_{(d_1,\dots,d_n)}(\tau)$ for the coefficients of the power series
\begin{equation}\label{hf}
\sum_{\tau=0}^\infty H_{(d_1,\dots,d_n)}(\tau)T^\tau=
\frac{\prod_{j=1}^n (1-T^{d_j})}{(1-T)^{n+1}}.
\end{equation}
It turns out that $H_{(d_1,\dots,d_n)}$  is the Hilbert function
of $\KK[x_0,x_1,\dots,x_n] / J$ when $J$ is an ideal  generated
by a regular sequence of $n$ homogeneous polynomials of degrees
$d_1,\dots, d_n$, that is, $H_{(d_1,\dots,d_n)}(\tau)$ is the
dimension as a $\KK$-vector space of the piece of degree $\tau$ in
$\KK[x_0,x_1,\dots,x_n] / J$; see \citet{Mac,Cha1}.

\begin{remark}\label{r2}
{}From the right-hand side of Identity (\ref{hf}), it is  easy to check that $H_{(d_1,\dots,d_n)}(\tau)<{\bf d}$
if $\tau<\rho,$ and $H_{(d_1,\dots,d_n)}(\tau)={\bf d}$
if $\tau\geq\rho.$
\end{remark}
\smallskip
If ${\rm Res}_{d_1,\dots,d_n}(f_{1d_1},\dots,f_{nd_n})\neq0$ holds,
Proposition \ref{affine} implies that  the
family of polynomials $f^0_1,\dots,f^0_n,x_0$ has no common
roots in projective space and so, ${\rm Res}_{d_1,\dots,d_n,1}(f^0_1,\dots,f^0_n,x_0)\neq0$.
But this implies that $f^0_1,\dots,f^0_n,x_0$ is a regular sequence in $\KK[x_0,\dots,x_n]$ and, in particular, $f^0_1,\dots,f^0_n$ is also a regular sequence in that ring.
Therefore,
$\dim \AA_{0\tau} = H_{(d_1,\dots,d_n)}(\tau)$.
\par
The next proposition shows a relationship between a monomial basis of the affine ring $\AA$ and bases of certain graded parts of the ring $\AA_0.$
This will allow us to state the condition for an arbitrary set $\MM$ to be a basis of $\AA$.
\smallskip
\begin{proposition}\label{affh}
If ${\rm Res}_{d_1,\dots,d_n}(f_{1d_1},\dots,f_{nd_n})\neq0,$ then
the following conditions are equivalent:
\begin{enumerate}
\item $\MM$ is a basis of $\AA$ as a ${\KK}$-vector space.
\item There exists $t_0\geq\max\{\delta,\rho\}$ such that
$\MM_{t_0}$ is a basis of $\AA_{0t_0}$ as a ${\KK}$-vector space.
\item For every
$t\geq\max\{\delta,\rho\},\,\MM_t$ is a basis of $\AA_{0t}$ as a ${\KK}$-vector space.
\end{enumerate}
\end{proposition}
\smallskip
\begin{remark}
We will see in Corollary \ref{bound} that a necessary condition for $\MM$ to be a basis of $\AA$ is that
$\delta\ge\rho$. Therefore, in the statement of Proposition \ref{affh} we can replace $\max\{ \delta,\rho\}$ with
$ \delta$.
\end{remark}
\smallskip

Now we will prove Proposition \ref{affh}.

\begin{proof}
Recall that the assumption ${\rm Res}_{d_1,\dots,d_n}(f_{1d_1},\dots,f_{nd_n})\neq0$
implies that $f^0_1,\dots, f^0_n$ is a regular sequence in
$\KK[x_0,\dots,x_n]$.

\smallskip

\par\noindent $\boxed{(1)\implies(3)}$ Let $t\geq\max\{\delta,\rho\}$ and consider a linear
combination of vectors in $\MM_t$ which lies in the ideal $( f^0_1,\dots,f^0_n):$
\begin{equation}\label{homm}
\sum_{i=1}^{\bf d} \lambda_i m_ix_0^{t-\deg(m_i)}=\sum_{j=1}^n A_j(x_0,\dots,x_n)f^0_j.
\end{equation}
Setting $x_0=1$ we get a linear combination of elements in $\MM$
which lies in $I.$ So, if $\MM$ is linearly independent, we get
that $\MM_t$ is linearly independent. As $t\geq\rho$ and
$f^0_1,\dots,f^0_n$ is a regular sequence, the dimension of
$\AA_{0t}$ is ${\bf d}$ and therefore, we conclude that $\MM_t$
is a basis of $\AA_{0t}$.
\par\noindent $\boxed{(3)\implies(1)}$
Consider a linear combination of $\MM$ as follows:
$$\sum_{i=1}^{\bf d} \lambda_i m_i=\sum_{j=1}^n a_j(x_1,\dots,x_n)f_j.$$
Let $t_0:=\max\{\delta,\,\rho,\,\deg(a_jf_j),\,j=1,\dots,n\}.$ Homogenizing
the linear combination up to degree $t_0,$ we have an equality like (\ref{homm}) with
$t_0$ instead of $t.$ As $\MM_{t_0}$ is
linearly independent, it turns out that $\lambda_i=0$ for $i=1,\dots,{\bf d}$. Then, $\MM$ is a linearly independent set.
Taking into account that $\dim(\AA) = {\bf d}$ it follows that it is a basis of $\AA.$
\par\noindent $\boxed{(3)\implies(2)}$ Obvious.
\par\noindent $\boxed{(2)\implies(3)}$ Consider the following exact complex of vector spaces:
$$0\to\ker{\phi_t}\to\AA_{0t}\stackrel{\phi_t}{\longrightarrow}\AA_{0(t+1)}
\to \left({\KK}[x_0,\dots,x_n] / (x_0,f^0_1,\dots,f^0_n)\right)_{t+1} \to0,$$
where $\phi_t(m)=x_0.m.$
As ${\rm Res}_{1,d_1,\dots,d_n}(x_0,f^0_1,\dots,f^0_n)\neq0,$ it turns
out that $\left(\KK[x_0,\dots,x_n] / (x_0,f^0_1,\dots,f^0_n)\right)_{t+1}=0$ if
$t\geq\rho.$ In addition, for $t\geq\rho$, we have that $\dim(\AA_{0t})=\dim(\AA_{0(t+1)}).$
So, $\phi_t$ is an isomorphism if $t\geq\max\{\rho,\delta\},$ and furthermore,
$\phi_t(\MM_t)=\MM_{t+1}.$ Then, $\MM_{t_0}$ is a basis of $\AA_{0t_0}$ for some $t_0\ge \max\{\delta, \rho\}$ if and only if $\MM_t$ is a basis of $\AA_{0t}$
for every $t\ge \max\{\delta, \rho\}$.
\end{proof}

The following result, which follows immediately from the proof of Proposition \ref{affh},
gives us a lower bound of the maximal degree one may expect from a monomial basis of $\AA.$

\begin{corollary}\label{bound}
If $\MM$ is a basis of $\AA,$ then $\delta(\MM)\geq\rho.$
\end{corollary}

\begin{proof}
Let $t<\rho,$ and suppose that $\MM$ is a basis of $\AA$ with $\delta=t.$
Proceeding as in the proof of $(1)\implies(3)$ in Proposition \ref{affh}, it follows that $\MM_t$ is linearly independent in $\AA_{0t}.$
But, from Remark \ref{r2}, we have that $\dim(\AA_{0t})<{\bf d}$ if $t<\rho,$ which is a contradiction.
\end{proof}
\smallskip
\begin{example}\label{example}
Let $f_1,f_2,f_3$ be generic polynomials of degree two in $\KK[x_1,x_2,x_3].$ In this
case, ${\bf d}= 2.2.2=8.$ It is well-known that
$$\MM:=\{1,x_1,x_2,x_3,x_1x_2,x_1x_3,x_2x_3,x_1x_2x_3\}$$ is generically a basis of $\AA$
(see for instance \citet{Mac}).
Observe that $\delta=3=\rho$ in this case. On the other hand,
Corollary \ref{bound} implies that there are no eight monomials
linearly independent in the set
$$\{1,x_1,x_2,x_3,x_1^2,x_2^2,x_3^2,x_1x_2,x_1x_3,x_2x_3\}.$$
This can be explained as follows: As $f^0_1,f^0_2,f^0_3$ is a
regular sequence, they must be linearly independent. So, the
dimension of the ${\KK}$-vector space they generate is $3$ and
hence, the dimension of $\AA_{02}$ is $10-3=7.$
\end{example}

\section{Subresultants by Means of Koszul Complexes}\label{pretres}

In this section we recall the theory of multivariate subresultants
for homogeneous polynomials as formulated in \citet{Cha1}; see
also \citet{Dem}.

First, we are going to introduce the crucial notion involved in
the definition of subresultants.

\subsection{The Determinant of an Exact Complex of Vector Spaces}

Let $K$ be a field and let $\mathbf{C}$ be an exact complex of
finitely generated $K$-vector spaces $F_i = K^{B_i}$, with bases
$B_i$, of the form
$$\mathbf{C} : 0 \rightarrow F_n \stackrel{\partial_n}{\rightarrow} F_{n-1}
\stackrel{\partial_{n-1}}{\rightarrow} \cdots
\stackrel{\partial_2}{\rightarrow} F_1
\stackrel{\partial_1}{\rightarrow} F_0 \to 0.$$ Then, there
exists a decomposition of the $K$-vector spaces $F_i$ which
enables us to associate with the complex $\mathbf{C}$ an element
$\Delta\in K$. This element $\Delta$ is called the {\em
determinant } of the complex  \citep[see][Appendix A]{GKZ}. In
order to obtain the decomposition, we can proceed as in
\citet{Dem,Cha1,GKZ}:

\medskip

\textsc{Ascending Decomposition}
\begin{itemize}
\item Set $I_1:= B_0$ and $V_1:=K^{I_1}$.
\item Since $\partial_1$ is onto, there exists a non-zero maximal minor
of the matrix of $\partial_1$. Choose such a non-zero minor, and
set $I_1'$ for the subset of $B_1$ corresponding to the elements
indexing the columns of the chosen submatrix and $I_2:= B_1 -
I_1'$. Then, if $V_1':= K^{I_1'}$ and $V_2:= K^{I_2}$, we have
$F_1 = V_2 \oplus V_1'$, and $\partial_1|_{V_1'} : V_1' \to V_1$
is an isomorphism.
\item For $i\ge 2$, consider $\partial_i^*:= \pi_{i-1} \circ \partial_i : F_i \to V_{i}$,
where $\pi_{i-1}$ is the projection from $F_{i-1}$ to $V_{i}$.
The map $\partial_i^*$ is onto, due to the exactness of
$\mathbf{C}$ and the chosen decomposition of $F_{i-1}$. Then, we
can choose a non-zero maximal minor of the matrix of
$\partial_i^*$ and consider the subset $I_{i}'$ of $B_i$ indexing
the columns of the chosen submatrix and $I_{i+1}:= B_i - I_{i}'$.
Setting $V_{i}':= K^{I_{i}'}$ and $V_{i+1}:= K^{I_{i+1}}$ we
obtain a decomposition $F_i = V_{i+1} \oplus V_{i}'$ such that the
restriction $\partial_i^*|_{V_{i}'}: V_{i}' \to V_{i}$ is an
isomorphism.
\item In the last step, we obtain a square matrix for $\partial_n^*$,
due to the fact that $\sum_{i=0}^n \dim(F_i)= 0$.
\end{itemize}

For every $1\le i \le n$, let $\phi_i:=
\partial_{i}^*|_{V_{i}'}: V_{i}' \to V_{i}$. The {\em
determinant} of the complex $\mathbf{C}$ (relative to the bases
$B_i$) is defined to be
$$\Delta :=\prod_{i=0}^{n-1} \det(\phi_{i+1})^{(-1)^i}.$$
We remark that $\Delta$ is (up to a sign) independent of the
choices made to perform the decomposition.

A second procedure to obtain a decomposition of a complex which
also enables us to compute its determinant, is the following:

\medskip

\textsc{Descending Decomposition}
\begin{itemize}
\item Set $I_n:= B_n$ and $V_n:= K^{I_n}$.
\item Since $\partial_n$ is into, there exists a non-zero maximal
minor of the matrix of $\partial_n$. Choose such a minor and
define $I_{n-1}\subset B_{n-1}$ to be the subset of elements of
$B_{n-1}$ indexing the rows not involved in this minor and
$I_n':= B_{n-1} - I_{n-1}$. Then we have a decomposition $F_{n-1}
= V_{n}' \oplus V_{n-1}$, where $V_n':= K^{I_n'}$ and $V_{n-1}:=
K^{I_{n-1}}$.

\item Note that, for $i\ge 1$, the previous construction for $i-1$
implies that ${\rm Im}(\partial_{n-i+1}) \cap V_{n-i} = 0$, and
therefore ${\rm Ker}(\partial_{n-i}) \cap V_{n-i} = 0$, that is,
the restriction of $\partial_{n-i}$ to $V_{n-i}$ is into. Then we
can iterate the process and choose a maximal non-zero minor of the
matrix of $\partial_{n-i}|_{V_{n-i}}$, and define $I_{n-i}'$ to
be the subset of $B_{n-i-1}$ indexing the rows of the chosen
submatrix and $I_{n-i-1}$ to be its complement in $B_{n-i-1}$. We
obtain a decomposition $F_{n-i-1}:= V_{n-i}' \oplus V_{n-i-1}$,
where $V_{n-i}':= K^{I_{n-i}'}$ and $V_{n-i-1}:= K^{I_{n-i-1}}$.

\item In the last step a square matrix is obtained, due to the
exactness of the complex.
\end{itemize}

As before, for every $1\le i \le n$, we define
$\phi_i:=\partial_{i}^*|_{V_{i}}: V_{i} \to V_{i}'$. It turns out
that \citep{GKZ,Cha1} the determinant of $\mathbf{C}$ relative to
the bases $B_i$ can also be computed as
$$\Delta :=\prod_{i=0}^{n-1} \det(\phi_{i+1})^{(-1)^i}.$$

\subsection{Subresultants}
Multivariate subresultants are defined as determinants of
generically exact Koszul complexes. Let $s\leq n+1$ and let
$P_1,\dots,P_s$ be generic homogeneous polynomials in $n+1$
variables $x_0,\ldots,x_n$ of respective degrees $d_1,\ldots,d_s$:
$$P_i(x_0,\ldots,x_n):=\sum_{|\alpha|=d_i}c_{i,\alpha}x^\alpha, \quad i=1,\ldots,s,$$
where the $c_{i,\alpha}$'s are new variables.
\par
In this case, $K$ is the field of fractions of
$A:=\ZZ\left[c_{i,\alpha},|\alpha|=d_i,\,i=1,\ldots,s\right]$. Set
$R:= A[x_0,x_1,\dots, x_n]$.
\par Let $\mathfrak{M}_t$ be the set of all monomials of degree $t$ in the variables
$x_0,\dots, x_n$, and let $S$ be a family of
$H_{d_1,\dots,d_s}(t)$ monomials in $\mathfrak{M}_t$. With this
data we can construct a complex $\mathbf{C}=\mathbf{C}_t^s$ which
is obtained by modifying the degree $t$ part of the Koszul
complex associated with $P_1,\dots, P_s$ as follows:
$$0 \rightarrow (\wedge^s R^s)_t \stackrel{\partial_s}{\rightarrow}
(\wedge^{s-1} R^{s})_t \stackrel{\partial_{s-1}}{\rightarrow}
\cdots \stackrel{\partial_2}{\rightarrow} (\wedge^1 R^s)_{t}
\stackrel{\varphi}{\rightarrow}A\langle \mathfrak{M}_t\setminus
S\rangle \rightarrow 0$$ equipped with the bases $B_k:=
\bigcup_{1\le i_1<\cdots<i_k\le s} \bigcup_{X^\alpha \in
\mathbb{M}_{t - d_{i_1} - \cdots - d_{i_k}}} X^\alpha e_{i_1}
\wedge \dots \wedge e_{i_k}.$

If this complex is generically exact (i.e. $\mathbf{C}\otimes K$
is exact as a complex of $K$-vector spaces), then the {\em
subresultant of $S$ with respect to the polynomials $P_1,\dots,
P_s,$} which will be denoted with $\Delta^t_S,$ is defined to be
the determinant of $\mathbf{C}\otimes K$ with respect to the
monomial bases; otherwise we set $\Delta^t_S:=0.$ As we have
$H_i(\mathbf{C}_t^s)=0$ for $i>0$ \citep{Jou,Cha1}, it turns out
that $\Delta^t_S$ is a {\em polynomial} in the coefficients of
the $P_i$'s which satisfies the following property \citep[Theorem
2]{Cha1}: Let $\mathbf{k}$ be any field, and
$\tilde{P_i}\in\mathbf{k}[x_0,\dots,x_n]_{d_i},\,i=1,\ldots,s.$
Then
$$\Delta^t_S(\tilde{P_1},\ldots,\tilde{P_s}) \ne 0
\iff J_t +\mathbf{k}\langle S\rangle =
\mathbf{k}[x_0,\dots,x_n]_t,$$ where $J_t$ is the degree $t$ part
of the ideal generated by the $\tilde{P_i}$'s.

\section{Monomial Bases and Subresultants}\label{tres}
In this section, we will relate our problem with multivariate subresultants.
\par
We set $s=n,$ and let $P_1,\dots, P_n$ be the homogeneous polynomials
$f^0_1,\dots,f^0_n$ defined above.
The following may be regarded as the main result of this section.
\begin{theorem}\label{mt}
Let $\MM\subset\KK[x_1,\dots,x_n]$ be a set of ${\bf d}$ monomials, and set
$t:=\delta(M).$ Let $\Delta^t_{\MM_t}$ be the subresultant of $\MM_t$ with respect to
$f^0_1,\dots,f^0_n.$ Then, $\MM$ is a basis
of $\AA$ if and only if
\begin{equation}\label{mbc}
P_{\MM,d_1,\dots,d_n}:={\rm
Res}_{d_1,\dots,d_n}(f_{1d_1},\dots,f_{nd_n})\Delta^t_{\MM_t}\neq0.
\end{equation}
\end{theorem}

\begin{proof}
If $\MM$ is a basis of $\AA$, the family $f_1,\dots,f_n$ has all its zeros in $\overline{\KK}^n,$ and therefore,
${\rm Res}_{d_1,\dots,d_n}(f_{1d_1},\dots,f_{nd_n})\neq0.$ In addition,
from Corollary \ref{bound} and Proposition \ref{affh} it follows that $\MM_t$ is a basis of $\AA_{0t},$
which implies that $\Delta^t_{\MM_t}\neq 0.$
\par
In order to prove the converse, we can apply Proposition \ref{affh}, as
${\rm Res}_{d_1,\dots,d_n}(f_{1d_1},\dots,f_{nd_n})\neq 0.$
The condition $\Delta^t_{\MM_t}\neq0$ implies that $\MM_t$ is a basis of $\AA_{0t}$ and then,
we conclude that
$\MM$ is a basis of $\AA$.
\end{proof}
\smallskip
\begin{example}
For $i=1,2,3,$ let $f_i:=\sum_{|\alpha|\leq2}c_{i,\alpha}x^\alpha$  be generic polynomials of degree two in $\KK[x_1,x_2, x_3]$, and let $\MM$ be as in example \ref{example}.
The subresultant $\Delta^3_{\MM_3}$ can be computed as the product of the determinants of the following two matrices:
$$\begin{pmatrix}
c_{1, 2, 0, 0} & c_{1, 0, 2, 0} & c_{1, 0, 0, 2}\\
c_{2, 2, 0, 0} & c_{2, 0, 2, 0} & c_{2, 0, 0, 2}\\
c_{3, 2, 0, 0} & c_{3, 0, 2, 0} & c_{3, 0, 0, 2}
\end{pmatrix}$$
and
{\footnotesize$$\begin{pmatrix}
c_{1,2,0,0} &0&0& c_{1,1,1,0} & c_{1,1,0,1} &0& c_{1,0,0,2} &0& c_{1,0,1,1} \\
0& c_{1,0,2,0} & 0 & c_{1,2,0,0} & 0 & c_{1,0,1,1} & 0 & c_{1,0,0,2} & c_{1,1,1,0} \\
0 & 0 & c_{1,0,0,2} &0& c_{1,2,0,0} & c_{1,0,2,0} & c_{1,1,0,1} & c_{1,0,1,1} & 0 \\
c_{2,2,0,0} & 0 & 0 & c_{2,1,1,0} & c_{2,1,0,1} & 0 & c_{2,0,0,2} & 0 & c_{2,0,1,1} \\
0 & c_{2,0,2,0} & 0 & c_{2,2,0,0} & 0 & c_{2,0,1,1} & 0 & c_{2,0,0,2} & c_{2,1,1,0} \\
0 & 0 & c_{2,0,0,2} &0& c_{2,2,0,0} & c_{2,0,2,0} & c_{2,1,0,1} & c_{2,0,1,1} & 0 \\
c_{3,2,0,0} & 0 & 0 & c_{3,1,1,0} & c_{3,1,0,1} & 0 & c_{3,0,0,2} & 0 & c_{3,0,1,1} \\
0 & c_{3,0,2,0} & 0 & c_{3,2,0,0} & 0 & c_{3,0,1,1} & 0 & c_{3,0,0,2} & c_{3,1,1,0} \\
0 & 0 & c_{3,0,0,2} & 0 & c_{3,2,0,0} & c_{3,0,2,0} & c_{3,1,0,1} & c_{3,0,1,1} & 0
\end{pmatrix}.$$}
\par
For a proof of this fact, see Theorem \ref{dependency} below.
\end{example}

\section{Factorization of Subresultants}\label{cuatro}
For several sets $\MM,$ the polynomial $P_{\MM,d_1,\dots,d_n}$
defined in (\ref{mbc}) depends only on the coefficients of
$f_{1d_1},\dots,f_{nd_n}$ and factorizes as a product of more
than two terms. For instance, \citet{Mac} showed that one can
decide whether
\begin{equation}\label{m0}
\MM^0:=\{x_1^{\alpha_1}\dots x_n^{\alpha_n},\,0\leq \alpha_i\leq d_i-1\}
\end{equation}
is a basis of $\AA$ by applying linear algebra on the
coefficients of the highest terms of $f_1, \dots, f_n$ \citep[see
also][{}]{BU}. The same has been done by \citet{BU} with
\begin{equation}\label{m1}
\MM^1:=\{x_1^{\alpha_1}x_2^{\alpha_2},\,0\leq \alpha_1< d_1,\,0\leq \alpha_2\leq d_1+d_2-2\alpha_1-2\},
\end{equation}
and with
$$
\{x_1^{\alpha_1}x_2^{\alpha_2}x_3^{\alpha_3},\,0\leq \alpha_1< d_1,\,0\leq \alpha_2< \min\left(
d_1,d_2,2(d_1-\alpha_1)-1\right),\qquad$$\vspace{-7mm}
$$\hspace{4.4cm}\,0\leq \alpha_3<d_1+d_2+d_3-2(\alpha_1+\alpha_2+1)\},
$$
for $n=2$ and $n=3$ respectively. This is not always the case, as the
following cautionary example shows.
\begin{example}
Consider $n=3.$ Set $d_1=d_2=d_3=2$ and write
$f_i:=\sum_{|\alpha|\leq2}c_{i,\alpha}\, x^\alpha$ for $i=1,2,3.$
Take $$\MM:=\{x_1^3,x_1,x_2,x_3,x_1x_2,x_1x_3,x_2x_3,x_1x_2x_3\}.$$
Then, $\Delta^3_{\MM_3}$ is the determinant of the following matrix:
{\scriptsize $$\left(\begin{array}{cccccccccccc}
c_{1, 0, 0, 0} & 0 & 0 & 0 & c_{2, 0, 0, 0} & 0 & 0 & 0 &
c_{3, 0, 0, 0} & 0 & 0 & 0\\
0 & 0 & c_{1, 0, 2, 0} & 0 & 0 & 0 & c_{2, 0, 2, 0} & 0 & 0 &
0 & c_{3, 0, 2, 0} & 0\\
0 & 0 & 0 & c_{1, 0, 0, 2} & 0 & 0 & 0 & c_{2, 0, 0, 2} & 0 &
0 & 0 & c_{3, 0, 0, 2}\\
c_{1, 2, 0, 0} & c_{1, 1, 0, 0} & 0 & 0 & c_{2, 2, 0, 0} &
c_{2, 1, 0, 0} & 0 & 0 & c_{3, 2, 0, 0} & c_{3, 1, 0, 0} & 0 & 0\\
c_{1, 0, 2, 0} & 0 & c_{1, 0, 1, 0} & 0 & c_{2, 0, 2, 0} & 0 &
c_{2, 0, 1, 0} & 0 & c_{3, 0, 2, 0} & 0 & c_{3, 0, 1, 0} & 0\\
c_{1, 0, 0, 2} & 0 & 0 & c_{1, 0, 0, 1} & c_{2, 0, 0, 2} & 0 &
0 & c_{2, 0, 0, 1} & c_{3, 0, 0, 2} & 0 & 0 & c_{3, 0, 0, 1}\\
0 & c_{1, 1, 1, 0} & c_{1, 2, 0, 0} & 0 & 0 & c_{2, 1, 1, 0} &
c_{2, 2, 0, 0} & 0 & 0 & c_{3, 1, 1, 0} & c_{3, 2, 0, 0} & 0\\
0 & c_{1, 1, 0, 1} & 0 & c_{1, 2, 0, 0} & 0 & c_{2, 1, 0, 1} &
0 & c_{2, 2, 0, 0} & 0 & c_{3, 1, 0, 1} & 0 & c_{3, 2, 0, 0}\\
0 & c_{1, 0, 2, 0} & c_{1, 1, 1, 0} & 0 & 0 & c_{2, 0, 2, 0} &
c_{2, 1, 1, 0} & 0 & 0 & c_{3, 0, 2, 0} & c_{3, 1, 1, 0} & 0\\
0 & c_{1, 0, 0, 2} & 0 & c_{1, 1, 0, 1} & 0 & c_{2, 0, 0, 2} &
0 & c_{2, 1, 0, 1} & 0 & c_{3, 0, 0, 2} & 0 & c_{3, 1, 0, 1}\\
0 & 0 & c_{1, 0, 0, 2} & c_{1, 0, 1, 1} & 0 & 0 &
c_{2, 0, 0, 2} & c_{2, 0, 1, 1} & 0 & 0 & c_{3, 0, 0, 2} &
c_{3, 0, 1, 1}\\
0 & 0 & c_{1, 0, 1, 1} & c_{1, 0, 2, 0} & 0 & 0 &
c_{2, 0, 1, 1} & c_{2, 0, 2, 0} & 0 & 0 & c_{3, 0, 1, 1} &
c_{3, 0, 2, 0}
\end{array}\right).$$}

With the aid of {\tt Maple} we have computed this determinant, which is an irreducible
polynomial depending on all the variables $c_{i,\alpha}.$
\end{example}

\bigskip

Set
\begin{equation}\label{ring}
\sum_{\tau=0}^\infty h_{(d_1,\dots,d_n)}(\tau)T^\tau=
 \frac{\prod_{j=1}^n (1-T^{d_j})}{(1-T)^{n}}.
\end{equation}
It turns out that $h_{d_1,\dots,d_n}$ is
the Hilbert function of the ideal generated by a regular sequence of $n$ homogeneous polynomials
in $n$ variables of degrees $d_1,\dots,d_n$ respectively.
\par
The following is the main result of this section:
\begin{theorem}\label{dependency}
Let $P_{\MM,d_1,\dots,d_n}$ be the polynomial defined in
(\ref{mbc}). Then, if $P_{\MM,d_1,\dots,d_n}$ is not identically
zero, the following conditions are equivalent:
\begin{itemize}
\item $P_{\MM,d_1,\dots,d_n}$ depends only on the coefficients of
$f_{1d_1},\dots,f_{nd_n}.$
\item For every $t=0,1,\dots,\rho,$ the cardinality of
$\MM\cap\KK[x_1,\dots,x_n]_t$ equals $h_{(d_1,\dots,d_n)}(t).$
\end{itemize}
If any of the above conditions hold, we have the following factorization:
\begin{equation}\label{factor}
\Delta^\delta_{\MM_\delta}=\prod_{t=\min\{d_i\}}^\rho D^t_{\MM\cap\KK[x_1,\dots,x_n]_t},
\end{equation}
where $D^t_S$ denotes the subresultant in $n$ variables of $S$ with respect to $f_{1d_1},\dots,f_{nd_n}.$
\end{theorem}
\begin{proof}
If $P_{\MM,d_1,\dots,d_n}$ depends only on the coefficients of
$f_{1d_1},\dots,f_{nd_n},$ we can set to zero all the
coefficients of $f_1,\dots, f_n$ not appearing in these leading
forms and work with this family of homogeneous polynomials
instead of $f_1,\dots,f_n.$ As $P_{\MM,d_1,\dots,d_n}$ is not
identically zero, we have that $\Delta^\delta_{\MM_\delta}$ is not
identically zero either and this implies that $\MM$ is a basis of
the homogeneous quotient ring $\KK[x_1,\dots,x_n] /
(f_{1d_1},\dots,f_{nd_n}).$ As the family
$f_{1d_1},\dots,f_{nd_n}$ is a regular sequence in
$\KK[x_1,\dots,x_n],$ it turns out that
 $\#\left(\MM\cap\KK[x_1,\dots,x_n]_t\right)=h_{(d_1,\dots,d_n)}(t)$ for any $t=0,\dots,\rho,$ and we are done.
\par
In order to prove the other implication, we will work with generic homogeneous polynomials. For each $i=1,\dots,n$ and $\alpha\in\NN_0^n$ with
$|\alpha|\leq d_i,$ introduce a variable $c_{i, \alpha}.$ Set
\begin{equation}\label{generique}
f_i(x_1,\dots,x_n):=\sum_{|\alpha|\leq d_i}c_{i,\alpha}\,x^{\alpha},\ i=1,\dots,n.
\end{equation}
We shall work in the field $\KK:=\QQ(c_{i, \alpha}).$ In this
situation we have that ${\rm
Res}_{d_1,\dots,d_n}(f_{1d_1},\dots,f_{nd_n})\neq0$ (see for
instance \citet{CLO2}) and, due to the universal property of
subresultants \citep{Cha1}, if $P_{\MM,d_1,\dots,d_n}\neq 0$ for
a given family of polynomials in \textit{any} field, then it will
not be zero for the generic family (\ref{generique}).

As before, set $f_i^0$ for the homogenization of the polynomial $f_i$ in $\KK[x_0,\dots,x_n].$ Consider the following $\KK$-linear map:
\begin{equation}\label{fhi}
\begin{array}{cccc}
\phi^\rho:&S^1_{\rho-d_1}\oplus\dots\oplus S^n_{\rho-d_n}&\to&S_{\rho} \\
&(p_1,\dots,p_n)&\mapsto&\sum_{i=1}^n p_i\,f^0_i,
\end{array}
\end{equation}
where $S_\rho:=\KK[x_0,x_1,\dots,x_n]_\rho,$ and for each $i=1,\dots,n,$
$$S^i_{\rho-d_i}:=\langle x_0^{\alpha_0}\dots x_n^{\alpha_n},\, \sum_{j=0}^n\alpha_j=\rho-d_i,\ \alpha_1<d_1,\dots,\alpha_{i-1}<d_{i-1}\rangle.$$
Let $M$ be the matrix obtained from the matrix of $\phi^\rho$ in
the monomial bases by  deleting the columns\footnote{As in
\citet{Mac}, the rows of $M$ are indexed by the monomial basis of
the domain.} indexed by the points in $\MM$ and let $M'$ be the
matrix obtained in the same way but using the set
\begin{equation}\label{tito}
S:=\{x_0^{\alpha_0}\dots x_n^{\alpha_n},\,|\alpha|=\rho,\, \alpha_i<d_i,\,i=1,\dots,n\}
\end{equation}
instead of $\MM$. It is well-known that $\det(M')\neq0$
\citep{Mac,Cha1}.
\par As the subresultant of $S$ with respect to $f_1^0, \dots, f_n^0$ is the determinant of $\mathbf{C}_t^S,$ it turns out that $\det(M')$ may be regarded
as a non-zero maximal minor in the last morphism of the complex
whose determinant is $\Delta^\rho_S.$
\par Starting with this maximal minor and using the ascending decomposition of the Koszul complex, it turns out that there exists an element
$\EE\in\KK$, which is actually
a polynomial in the $c_{i,\alpha},$ such that $\det(M')=\EE\,\Delta^\rho_S.$ As $\det(M')\neq0,$ then $\EE\neq0.$
\par
This $\EE$ is a product of complementary minors in
$\mathbf{C}_t^S.$ Starting now with these minors from the left
and applying the descending decomposition of the Koszul complex,
one can see that, as in \citet{Cha1},
$\det(M)=\EE\,\Delta^\rho_\MM,$ as the complex whose determinant
is $\Delta^\rho_\MM$ is the same as the one whose determinant is
$\Delta^\rho_S$ except in the last map.
\par

Set $\MM(t):=\MM\cap\KK[x_1,\dots,x_n]_t,\,t=0,1,\dots,\rho,$ and
suppose w.l.o.g. that $d_1\leq d_i,\,i=2,\dots,n.$ As
$\#\MM(t)=h_{d_1,\dots,d_n}(t),$ proceeding as in \citet{Mac}, it
follows that --ordering appropriately its rows and columns-- the
matrix $M$ has the following block structure:
\begin{equation}\label{tit}
\left(\begin{array}{cccccc}
M_\rho&*&*&*\\
0&M_{\rho-1}&*&*\\
0&0&\ddots &*\\
0&0&\dots &M_{d_1}
\end{array}\right),
\end{equation}
where $M_t$ is the square matrix obtained by deleting the columns indexed by the monomials in $\MM(t)$ in the matrix of the $\KK$-linear map:
$$
\begin{array}{cccc}
\phi_t:&S^{1*}_{t-d_1}\oplus\dots\oplus S^{n*}_{t-d_n}&\to&S^*_{t} \\
&(p_1,\dots,p_n)&\mapsto&\sum_{i=1}^n p_i\,f_{id_i}.
\end{array}
$$
Here $S^*_t:=\KK[x_1,\dots,x_n]_t,$ and for each $i=1,\dots,n,$
$$S^{i*}_{t-d_i}:=\langle x_1^{\alpha_1}\dots x_n^{\alpha_n},\, \sum_{k=1}^n\alpha_k=t-d_i,\ \alpha_1<d_1,\dots,\alpha_{i-1}<d_{i-1}\rangle.$$
Then, we have that $\det(M)=\prod_{t=d_1}^\rho\det(M_t),$ which
shows that $\det(M)$ depends only on the coefficients of
$f_{id_i},i=1,\dots,n.$ Furthermore, $\det(M_t)=\EE_t
D^t_{\MM\cap\KK[x_1,\dots,x_n]_t}$ for $t=0,\dots, \rho$, and the
extraneous factor $\EE$ has also a block structure compatible
with the one given in (\ref{tit}), that is,
$\EE=\prod_{t=d_1}^\rho\EE_t$; see \citet{Mac,Cha2}. This
completes the proof of the theorem.
\end{proof}

\medskip
\begin{corollary}
If $P_{\MM,d_1,\dots,d_n}$ is not identically zero and depends
only on the coefficients of $f_{1d_1},\dots,f_{nd_n},$ then
$\delta(\MM)=\rho.$
\end{corollary}

\section{Simple Roots and Generalized Vandermonde Determinants}\label{sr}\label{cinco}

In this section, we will study a result by \citet{Mac} concerning
the structure of a generalized Vandermonde determinant associated
with the monomial set $\MM^0$ and, with the aid of subresultants,
we will extend it to arbitrary sets of monomials with cardinality
$\mathbf{d}$. This will make apparent the relationship between the
non-vanishing of the generalized Vandermonde determinant
associated with a set of monomials $\MM$ and the fact that $\MM$
is a basis of the quotient algebra $\AA$ in the case of a
polynomial system with simple roots.

We will work in the generic field $\KK=\QQ(c_{i,\alpha}),$ and
with the family (\ref{generique}). Let
$V(f_1,\dots,f_n)=\{\xi_1,\dots,\xi_{\bf
d}\}\subset\overline{\KK}^n,$ and set $\MM^0=\{m_1,\dots,m_{\bf
d}\}$ (recall that $\MM^0$ was defined in (\ref{m0})). Let $M_0$
be the ${\bf d}\times{\bf d}$ matrix whose rows (resp. columns)
are indexed by the elements of $V(f_1,\dots,f_n)$ (resp. $\MM^0$),
such that the element indexed by $(\xi_i,m_j)$ is the evaluation
of $m_j$ at $\xi_i,$ that is, $M_0:= \left(m_j(\xi_i)\right)_{1\le
i, j\le n}$.
\par
In \cite[Section 10]{Mac}, it is proven that
\begin{equation}\label{disp}
{\det(M_0)}^2={\bf c}\,{\mathcal J}\frac{({\Delta^\rho_{\MM^0_\rho}})^2}{{{\rm Res}_{(d_1,\dots,d_n)}(f_{1d_1},\dots,f_{nd_n})}^{\rho+1}},
\end{equation}
where ${\mathcal J}:=\prod_{i=1}^{{\bf d}} J(\xi_i)$ (here
$J:=\det\left({\partial f_i}/ {\partial x_j}\right)_{1\leq i,j\leq
n}$ is the Jacobian of the sequence $f_1,\dots,f_n$), and ${\bf
c}\in \mathbb{Q}$ is a numerical constant depending only on $n$
and the degrees $d_1,\dots, d_n$.

\smallskip
The constant ${\bf c}$ in (\ref{disp}) has an explicit expression in terms of $d_1,\dots, d_n$:

\begin{lemma}\label{constante}
$$\,{\bf c}=(-1)^{E_n(d_1,\dots, d_n)},$$ where
$$E_n(d_1,\dots, d_n):=\sum_{j=1}^n d_1\dots d_{j-1} \frac{(d_j-1)d_j}{2}d_{j+1} \dots d_n.$$
\end{lemma}

\begin{proof}
First, observe that a system $f_1, \dots, f_n$ having the property that $f_{id_i}= x_i^{d_i}$ for
$i=1,\dots, n$,  verifies ${\rm
Res}_{(d_1,\dots,d_n)}(f_{1d_1},\dots,f_{nd_n})=1$ and
$({\Delta^\rho_{\MM^0_\rho}})^2=1$, as both polynomials depend only
on the coefficients of $f_{1d_1},\dots,f_{nd_n}$ (see Theorem
\ref{dependency} above). Therefore, the numerical factor ${\bf
c}$ can be obtained from identity (\ref{disp}) by specializing
the coefficients of $f_i$ in such a way that
$f_{id_i}=x_i^{d_i},\,i=1\dots,n.$ If this is the case, we get
\begin{equation}\label{formc}
{\bf c}= \dfrac{{\det(M_0)}^2}{\mathcal J}.
\end{equation}
The theorem will be proved by induction on $n$.
\par
First, we fix some notation. We denote by $c_n(d_1, \dots, d_n)$
the numerical factor associated with $n$ and degrees $d_1,\dots,
d_n$. If $f_1,\dots, f_n$ is a system of polynomials in $n$
variables of degrees $d_1,\dots, d_n,$ we denote by
$\mathcal{M}_n(f_1,\dots, f_n)$ the matrix $M_0$ associated with
the system $f_1, \dots, f_n$ and the set $\MM^0$, and we set
${\mathcal J}_n(f_1,\dots, f_n) :=\prod_{i=1}^{{\bf d}} J(\xi_i).$
\par
For $n=1$, set $d_1 = d$ for a positive integer and let $f_1:=
x_1^d -1 $. We have that $V(f_1) = \{ \xi_1, \dots, \xi_d\}$ is
the set of $d$th roots of unity. The matrix $M_0$ is the
Vandermonde matrix associated with the roots of $f_1$ and so,
$\det(M_0)^2 = {\rm disc}(f_1) = (-1)^{d-1+\frac{d(d-1)}{2}}
d^d$. In addition, ${\mathcal J} = (-1)^{d-1} d^d$. Then we
conclude from identity (\ref{formc}) that
\begin{equation*}
c_1(d) = (-1)^\frac{d(d-1)}{2}.
\end{equation*}

Assume now that the formula holds for systems of $n$ polynomials in $n$ variables and consider $n+1$ polynomials in $n+1$ variables.
\par
$\bullet$ For degrees $d_1,\dots, d_n, 1$: Set $f_i := x_i^{d_i} - 1$ for $i=1,\dots, n$, and $f_{n+1} := x_{n+1}$. We have
$$
V(f_1,\dots, f_{n+1})=\{ (\eta_1,\dots, \eta_n, 0): \eta_i^{d_i}=1,\, 1\le i \le n\},
$$ and so, it is straightforward to check that
\begin{eqnarray*}
\mathcal{M}_{n+1}(f_1,\dots, f_n, f_{n+1}) &=&  \mathcal{M}_{n}(x_1^{d_1}-1,\dots, x_n^{d_n}-1),\\
{\mathcal J}_{n+1}(f_1,\dots, f_n, f_{n+1}) &=&  {\mathcal J}_{n}(x_1^{d_1}-1,\dots, x_n^{d_n}-1).
\end{eqnarray*}
Identity (\ref{formc}) implies
$$
c_{n+1}(d_1,\dots, d_n, 1) = c_n(d_1,\dots, d_n),
$$
and the formula holds.
\par
$\bullet$ For degrees $d_1,\dots, d_n, d_{n+1}+1$: Set $f_i:=
x_i^{d_i} - 1$ for $1\le i \le n$, and $f_{n+1} :=
x_{n+1}^{d_{n+1}+1}- x_{n+1}$. Then, $V(f_1,\dots, f_{n+1}) = V_1
\cup V_2$, where $V_1 = V(x_1^{d_1} - 1,\dots, x_n^{d_n} - 1)
\times \{ 0\}$ and $V_2 = V(x_1^{d_1} - 1,\dots, x_n^{d_n} -
1)\times \{ \eta \in \overline{\KK} : \eta^{d_{n+1}} =1\}$.
Arranging the monomials in $\MM^0$ so that those which do not
depend on the variable $x_{n+1}$ come first and the roots of the
system so that those in $V_1$ come first, it follows that
$\mathcal{M}_{n+1}(f_1,\dots, f_{n+1}) $ has the following block
structure:
$$  \left(\begin{array}{cc}
\mathcal{M}_{n}(x_1^{d_1}-1,\dots, x_n^{d_n}-1) & 0 \\
* & \mathcal{M}_{n+1}'(x_1^{d_1}-1,\dots, x_n^{d_n}-1, x_{n+1}^{d_{n+1}}-1)
\end{array}\right)$$
where $\mathcal{M}_{n+1}'(x_1^{d_1}-1,\dots, x_n^{d_n}-1,
x_{n+1}^{d_{n+1}}-1)$ is a matrix differing from
$\mathcal{M}_{n+1}(x_1^{d_1}-1,\dots, x_n^{d_n}-1,
x_{n+1}^{d_{n+1}}-1)$ only in a factor by a $d_{n+1}$th root of
unity in each row. Moreover, each root of unity appears in
exactly $d_1\dots d_n$ rows. Taking into account that the product
of all the $d_{n+1}$th roots of unity equals $(-1)^{d_{n+1}-1}$,
it follows that $\big(\det \mathcal{M}_{n+1}(f_1,\dots,
f_{n+1})\big)^2$ equals the product
$$
\big(\det \mathcal{M}_{n}(x_1^{d_1}-1,\dots, x_n^{d_n}-1)\big)^2
 \big(\det \mathcal{M}_{n+1}(x_1^{d_1}-1,\dots, x_n^{d_n}-1,
x_{n+1}^{d_{n+1}} -1)\big)^2.
$$
On the other hand, the Jacobian of the polynomial system
$f_1,\dots, f_n, f_{n+1}$ is $J= d_1x_1^{d_1-1}\dots d_n
x_n^{d_n-1}((d_{n+1}+1) x_{n+1}^{d_{n+1}} -1)$ and then, for
every $\xi \in V_1$, $J(\xi) = (-1) J(x_1^{d_1}-1,\dots,
x_n^{d_n}-1)(\xi)$ and, for every $\xi \in V_2$, $J(\xi) =
\xi_{n+1} J(x_1^{d_1}-1,\dots, x_{n+1}^{d_{n+1}}-1)(\xi)$. Then,
it follows easily that
$$\prod_{\xi \in V_1} J(\xi) = (-1)^{d_1\dots d_n} {\mathcal J}_n(x_1^{d_1}-1,\dots, x_n^{d_n}-1),$$
$$\prod_{\xi \in V_2} J(\xi) = (-1)^{d_1\dots d_n (d_{n+1}-1)}
{\mathcal J}_{n+1}(x_1^{d_1}-1,\dots,
x_n^{d_n}-1,x_{n+1}^{d_{n+1}}-1 )$$ and so, ${\mathcal
J}_{n+1}(f_1,\dots, f_{n+1})$ equals
$$(-1)^{d_1\dots d_n d_{n+1}} {\mathcal J}_n(x_1^{d_1}-1,\dots, x_n^{d_n}-1) {\mathcal J}_{n+1}(x_1^{d_1}-1,\dots, x_n^{d_n}-1,x_{n+1}^{d_{n+1}}-1 ).$$
{}From the expressions for $\mathcal{M}_{n+1}$ and ${\mathcal
J}_{n+1}$, we deduce:
$$
c_{n+1}(d_1,\dots, d_n, d_{n+1}+1)= (-1)^{d_1\dots d_n d_{n+1}} c_n(d_1,\dots,d_n)c_{n+1}(d_1,\dots, d_n, d_{n+1}).
$$
Thus, the inductive assumption implies that $c_{n+1}(d_1,\dots, d_n, d_{n+1}+1)= \pm 1$.
More precisely, the exponent
${E_{n+1}(d_1,\dots, d_n, d_{n+1}+1)}$ giving the sign equals
$$
d_1\dots d_n d_{n+1} + E_n(d_1,\dots, d_n) +
E_{n+1}(d_1,\dots, d_n, d_{n+1}) ={}
$$
$$ {}=
\sum_{j=1}^{n+1} d_1\dots d_{j-1} \frac{(d_j-1)d_j}{2}\, d_{j+1}
\dots d_n d_{n+1}.$$
\end{proof}

\smallskip
Let $\MM$ be any set of monomials of cardinality ${\bf d},$ and let
$M:=M(\MM)$ be the matrix defined as $M_0$ but with the columns indexed by the elements of $\MM.$
The main result of this section is an expression similar to (\ref{disp}) for $M$:
\begin{theorem}\label{main}
$${\det(M(\MM))}^2= \pm \,{\mathcal J}\frac{{(\Delta^\delta_{\MM_\delta})}^2}{{{\rm Res}_{(d_1,\dots,d_n)}(f_{1d_1},\dots,f_{nd_n})}^{2\delta-\rho+1}}.$$
\end{theorem}

\medskip
The following result will be needed in the proof of Theorem \ref{main}.
\begin{lemma}\label{ff}
For any $t\geq\delta=\delta(\MM),$
$$\Delta^{t}_{\MM_{t}}=\Delta^{\delta}_{\MM_\delta}\,{\rm Res}_{(d_1,\dots,d_n)}(f_{1d_1},\dots,f_{nd_n})^{t-\delta}.$$
\end{lemma}

\begin{proof}
It is enough to prove the result for $t=\delta+1$ and $\delta\geq\rho$ (otherwise, both subresultants are identically zero and the claim holds).

Consider the morphisms for computing $\Delta^\delta_{\MM_\delta}$ and
$\Delta^{\delta+1}_{\MM_{\delta+1}}$ as in (\ref{fhi}):
\begin{equation}\label{sqc}
\begin{array}{ccc}
S^1_{\delta-d_1}\oplus\dots\oplus S^n_{\delta-d_n}&\stackrel{\phi^\delta}{\to}&S_{\delta} \\
\downarrow&&\downarrow \\
S^1_{\delta+1-d_1}\oplus\dots\oplus S^n_{\delta+1-d_n}&\stackrel{\phi^{\delta+1}}{\to}&S_{\delta+1},
\end{array}
\end{equation}
where the vertical maps are multiplication by $x_0.$ It is straightforward to check that the diagram (\ref{sqc}) commutes. For $i=\delta,\delta+1,$ let $M^i$ be the matrix
of $\phi^{i}$ where we have deleted the columns indexed by those $m\in\MM_{i}$. If we order the rows and columns of $M^{\delta+1}$
in such a way that the monomials having degree zero in $x_0$ come first, it is easy to see that this matrix has the following structure:
$$\left(\begin{array}{cc}
M_{\delta+1}&* \\
0 & M^{\delta}
\end{array}
\right),
$$
where $M_{\delta+1}$ has been defined in the proof of Theorem \ref{dependency}.
\par As $\delta+1>\rho,$ there exists a polynomial $\EE_1\in \QQ[c_{i,\alpha}]$ such that  $\det(M_{\delta+1})={\rm
Res}_{(d_1,\dots,d_n)}(f_{1d_1},\dots,f_{nd_n})\EE_1$
\citep{Mac}. Besides, there are also elements $\EE_2$ and $\EE$
such that $\det(M^\delta)=\Delta^\delta_{\MM_\delta}\EE_2$ and
$\det(M^{\delta+1})=\Delta^{\delta+1}_{\MM_{\delta+1}}\EE.$ As in
the proof of Theorem \ref{dependency}, we use the block structure
of the extraneous factor $\EE$ \citep{Mac,Cha2}, and it turns out
that $\EE=\EE_1\EE_2.$
\end{proof}

\bigskip
\begin{pot}
Let $\delta=\delta(\MM)$. If $\Delta^\delta_{\MM_\delta}=0,$ it
follows that the same holds for $\det(M(\MM))$.
\par
If this is not
the case, consider the following complex of
$\overline{\KK}$-vector spaces:
\begin{equation}\label{ufa}
0\to S^1_{\delta-d_1}\oplus\dots\oplus S^n_{\delta-d_n}\stackrel{\phi}{\to}S_{\delta}\stackrel{\psi}{\to}\overline{\KK}^{\bf d}\to0,\end{equation}
where $S_\delta:=\overline{\KK}[x_0,x_1,\dots,x_n]_\delta$ and, as before,
$$\begin{array}{l}
S^i_{\delta-d_i}:=\langle x_0^{\alpha_0}\dots x_n^{\alpha_n},\, \sum_{j=0}^n\alpha_j=\delta-d_i,\ \alpha_1<d_1,\dots,\alpha_{i-1}<d_{i-1}\rangle_{\overline{\KK}},\\[2mm]
\phi(p_1,\dots,p_n):=\sum_{i=1}^n p_if^0_i,\\[2mm]
\psi(p(x)):=(p(1,\xi_1),\dots,p(1,\xi_{\bf d})).
\end{array}$$
It is easy to see that the complex (\ref{ufa}) is exact.
If $\MM'$ is another set of ${\bf d}$ elements such that $\delta(\MM')\leq\delta(\MM)$ and $\det(M(\MM'))\neq0,$ we denote with $D(\MM'_\delta)$
(resp. $D(\MM_\delta)$) the determinant of
the matrix of $\phi$ in the monomial bases where we have deleted the columns indexed by those monomials lying in $\MM'_\delta$ (resp. $\MM_\delta$).
Then, considering the determinant of the complex (\ref{ufa}), we have the following:
$$\frac{D(\MM_\delta)}{\det(M(\MM))}=\pm\frac{D(\MM'_\delta)}{\det(M(\MM'))}.$$
As in the proof of Theorem \ref{dependency}, it turns out that
$D(\MM'_\delta)=\EE\,\Delta^\delta_{\MM'_\delta}$ and $D(\MM_\delta)=\EE\,\Delta^\delta_{\MM_\delta}$, with the same extraneous factor $\EE.$ Therefore
$$\frac{\Delta^\delta_{\MM_\delta}}{\det(M(\MM))}=\pm\frac{\Delta^\delta_{\MM'_\delta}}{\det(M(\MM'))}.$$
Taking as $\MM'$ the set $\MM^0$, it follows that
$$\left(\frac{\Delta^\delta_{\MM_\delta}}{\det(M(\MM))}\right)^2=
\left(\frac{\Delta^\delta_{\MM^0_\delta}}{\det(M_0)}\right)^2
=\left(\frac{\Delta^\rho_{\MM^0_\rho}{{\rm Res}_{(d_1,\dots,d_n)}(f_{1d_1},\dots,f_{nd_n})}^{\delta-\rho}}{\det(M_0)}\right)^2,$$
where the last equality holds for Lemma \ref{ff}.

Now, the claim is an immediate consequence of identity (\ref{disp}) and Lemma \ref{constante}.
\end{pot}

\section{An Overview of the B\'ezout Construction of the Resultant}\label{seis}

In this section we will compare several results obtained by
\citet{BU} with ours. This will allow us to clarify the B\'ezout
construction of the resultant.
\par
In \cite[Section 4]{BU}, the matrix $M_0$ defined at the
beginning of Section \ref{sr} is introduced (it is denoted as
$V$) and the structure of ${\det(M_0)}^2$ is studied. Following
\citet{Mac}, it is stated that
$$\det(M_0)^2 = \Upsilon \mathcal{J},$$
where $\mathcal{J}$ is as defined in Section \ref{cinco} of this
paper. Furthermore, it is claimed that $\Upsilon$ is a rational
function in the coefficients of the leading forms of the
polynomials $f_1,\dots,f_n$ whose numerator is a product of
$\rho$ polynomials in these coefficients.

In our notation, identity (\ref{disp}) and Lemma \ref{constante}
imply that
$$\Upsilon=\pm\frac{({\Delta^\rho_{\MM^0_\rho}})^2}{{{\rm
Res}_{(d_1,\dots,d_n)}(f_{1d_1},\dots,f_{nd_n})}^{\rho+1}}.$$
Moreover, the fact stated in \citet{BU} about the factorization
of the numerator of $\Upsilon$ is Theorem \ref{dependency} of the
present paper applied to $\MM^0$  \citep[see also][Section
10]{Mac}. Finally, let us observe that the irreducible factors of
the numerator and the denominator of $\Upsilon$ and of the
polynomial $P_{\MM^0,d_1,\dots,d_n}$ defined in Theorem \ref{mt}
are the same and, therefore, due to our main result we have that
$\Upsilon\neq0$ if and only if $\MM^0$ is a basis of $\AA.$
\par

Also, the structure of $\det(M(\MM^1))^2$ is studied in
\cite[Theorem 5.1]{BU} in the bivariate case (see the definition
of $\MM^1$ in (\ref{m1})). We point out a mistake in formula
$(5.30)$ of \citet{BU}, which is incorrect if the degrees of the
input polynomials are different. This follows straightforwardly
due to the fact that $\det(M(\MM^1))^2$ has degree zero in the
coefficients of $f_1,\dots,f_n,$ and if $n=2,$ then ${\mathcal
J}$ has degree $2d_1d_2$ in these coefficients and the $k$th
classical subresultant has degree
$d_1+d_2-2k,\,k=1,\dots,\min(d_1,d_2).$ If $d_1<d_2,$ it turns
out that the $k$th classical subresultant is the multivariate
subresultant of $\MM^1_{\rho-k+1}$ with respect to
$f_{1d_1},f_{2d_2}$ if $1\leq k \leq d_1-1$ \citep{Cha1}. It
remains to compute the multivariate subresultant of $\MM^1_t$ for
those degrees $t$ such that $d_1\leq t < d_2.$ This is easily
seen to be equal to $c_{1,(d_1,0)}^{t+1-d_1}.$ Hence, we have the
following
\begin{proposition}\label{correc1}
$$\Upsilon = {\bf
c}\frac{{(\rr_1\dots\rr_{d_1-1})}^2c_{1,(d_1,0)}^{(d_2-d_1)(d_2-d_1+1)}}{{{\rm
Res}_{(d_1,d_2)}(f_{1d_1},f_{2d_2})}^{\rho+1}},
$$
where $\rr_i$ is the classical $i$-subresultant and ${\bf c}$ is the constant of Lemma \ref{constante}.
\end{proposition}

Concerning the reducibility problem (that is, given a family of
polynomials $f_1,\dots, f_n$ with respective degrees $d_1,\dots,
d_n$ and a set of monomials $\MM$ with cardinality $\mathbf{d} =
d_1\dots d_n$, decide whether every polynomial is a linear
combination of $\MM$ when reduced modulo the ideal $(f_1,\dots,
f_n)$), in Section $5$ of \citet{BU}, a reduction algorithm with
respect to $\MM^0$ and $\MM^1$ is presented by solving a
succession of linear systems whose coefficients depend rationally
on the leading forms of the input polynomials. One can easily
check that the matrices of these linear systems can be regarded
as subresultant matrices. Indeed, in \cite[Theorem 5.1]{BU},
reduction modulo $\MM^1$ is completely characterized in terms of
the classical subresultants if $n=2.$

In \cite[Theorem $5.2$]{BU} it is claimed that, for three
polynomials of equal degree $d,$ it is sufficient for
reducibility that $2d-1$ determinants are non-zero. However, as a
result of Theorem \ref{dependency}, we get that $2d-2$ conditions
suffice. This can be verified following the approach by \citet{BU}
in detail: it turns out that the linear systems they consider have
determinants which are rational functions involving
subresultants, and that the condition arising in the last system
in their algorithm is redundant. Also, in \cite[Theorem 5.3]{BU}
it is shown that the first $d$ conditions of the $2d-1$ needed in
their reduction algorithm can be rewritten in terms of the nested
minors of the Macaulay matrix of the initial forms of the
polynomials. This follows straightforwardly in our framework, due
to the structure of the Macaulay matrix given in (\ref{tit}) and
the fact that, for $d\leq t\leq 2d-1,$
$\det(M_t)=D^t_{\MM\cap\KK[x_1,\dots,x_n]_t},$ i.e. there are no
extraneous factors \citep{Mac}.

Similar remarks can be made about the general approach they
present in \cite[Section $5.3.$]{BU}.
\par

Finally, we will answer negatively the Rank Conjecture posted in
\cite[Section 4]{BU}. Let $f_1,\dots,f_n$ be polynomials such that
$\MM^0$ is a basis of $\AA.$ Let $g\in\KK[x_1,\dots,x_n],$ and
let us denote with $\BB$ the matrix of the following linear map in
the basis $\MM^0$:
\begin{equation}\label{bezz}
\begin{array}{ccc}
\AA&\to&\AA \\
p(x)&\mapsto&p(x)\,g(x).
\end{array}
\end{equation}
It is a well-known fact \citep[see][{}]{CLO2,BU} that if $V(g)\cap
V(f_1,\dots,f_n)=\emptyset,$ then the determinant of $\BB$ equals
the dense resultant of the family $f_1,\dots,f_n,g$ up to a
constant. Suppose now that $V(g)\cap
V(f_1,\dots,f_n)=\{p_1,\dots,p_s\},$ and for each $i=1,\dots,s,$
we denote with $l_i$ the minimum between the multiplicity of
$p_i$ as a zero of $V(f_1,\dots,f_n)$ and the multiplicity of
$p_i$ as a zero of $g.$ The Rank Conjecture asserts that the rank
of $\BB$ should be equal to ${\bf d}-\sum_{i=1}^s l_i.$

\par This conjecture is
not true in general. For instance, we can take $f_1,\dots,f_n$
homogeneous polynomials of respective degrees $d_1,\dots,d_n$
such that the specialization of $P_{\MM^0\!,d_1,\dots,d_n}$ in
the coefficients of this family is not identically zero. This
implies that the only zero of the affine variety
$V(f_1,\dots,f_n)$ is the zero vector with multiplicity ${\bf
d}.$ Moreover, $\MM^0$ is a basis of $\AA,$ which is a graded
ring of finite dimension with $\AA_t=0$ for $t>\rho.$ Let $g$ be
any homogeneous polynomial of degree $d.$ According to the Rank
Conjecture, the kernel of $\BB$ should have dimension equal to
$\min\{{\bf d},d\},$ which is true if $d=0$ or $d>{\bf d},$ but
not in general. A straightforward computation shows that
$\AA_t\subset\ker(\BB)$ if $t>\rho-d,$ so
$$\dim\left(\ker(\BB)\right)\geq\sum_{j=\rho-d+1}^\rho h_{(d_1,\dots,d_n)}(j),$$ and this number may be greater than $d.$ For instance, if $d=2,\,d_i>3,$ we have that
$$h_{(d_1,\dots,d_n)}(\rho-1)+h_{(d_1,\dots,d_n)}(\rho)=n+1,$$
which is greater than $2$ unless $n=1.$
\par

\begin{ack}
We are grateful to P. Bikker and A.~Yu. Uteshev for providing us
updated versions of their joint work, to Laurent Bus\'e for
helpful comments on a preliminar version of this paper, and to
the anonymous referees for their useful suggestions.
\par The
first author was supported by the Miller Institute for Basic
Research in Science, in the form of a Miller Research Fellowship
(2002--2005). The second author was partially supported by
CONICET, grant PIP 2461 (2000-2002), and Universidad de Buenos
Aires, grant X198 (2001-2003).
\end{ack}


\begin{thebibliography}{ }

\bibitem[B\'ezout(1779)]{Bez}
B\'ezout, \'E., 1779.
\newblock{ Th\'eorie g\'en\'erale des \'equations alg\'ebriques.}
\newblock Ph.-D. Pierres: Paris.

\bibitem[Bikker \& Uteshev(1999)]{BU}
Bikker, P., Uteshev, A.~Yu., 1999.
\newblock{On the B\'ezout construction of the resultant.}
\newblock J. Symb. Comput. 28, No.1-2, 45-88.


\bibitem[Chardin(1994a)]{Cha2}
Chardin, M., 1994a.
\newblock{Formules \`a la Macaulay pour les sous-r\'esultants en plusieurs variables.}
\newblock C. R. Acad. Sci. Paris S\'er. I Math.  319,  no. 5, 433--436.

\bibitem[Chardin(1994b)]{Cha3}
Chardin, M., 1994b.
\newblock{Sur l'ind\'ependance lin\'eaire de certains mon\^omes modulo des polyn\^omes g\'en\'eriques.}
\newblock C. R. Acad. Sci. Paris S\'er. I Math.  319,  no. 10, 1033--1036.

\bibitem[Chardin(1995)]{Cha1}
Chardin, M., 1995.
\newblock{Multivariate subresultants.}
\newblock J. Pure Appl. Algebra 101, No.2, 129-138.

\bibitem[Cox et al.(1996)]{CLO1}
Cox, D., Little, J., O'Shea,  D., 1996.
\newblock{Ideals, varieties, and algorithms.
An introduction to computational algebraic geometry and commutative algebra.}
\newblock Undergraduate Texts in Mathematics. New York, NY: Springer.
xiii.

\bibitem[Cox et al.(1998)]{CLO2}
Cox, D., Little, J., O'Shea,  D., 1998.
\newblock{Using algebraic geometry}.
\newblock  Graduate Texts in Mathematics, 185. Springer-Verlag, New
York.

\bibitem[Demazure(1984)]{Dem}
Demazure, M., 1984.
\newblock{Une d\'efinition constructive du r\'esultant.}
\newblock Notes Informelles de Calcul Formel {\bf 2}, pr\'epublication du Centre de Math\'ematiques de l\'Ecole
Polytechnique.

\bibitem[Emiris \& Rege(1994)]{ER}
Emiris, I.~Z., Rege,  A., 1994.
\newblock{Monomial bases and polynomial system solving.}
\newblock In Proceedings of International Symposium on Symbolic
and Algebraic Computation, Oxford, 114--122. ACM: New York.

\bibitem[Gel'fand et al.(1994)]{GKZ}
Gel'fand, I.~M., Kapranov, M.~M., Zelevinsky, A.~V., 1994.
\newblock{Discriminants, resultants, and multidimensional determinants.}
\newblock Mathematics: Theory \& Applications. Birkh\"auser Boston, Inc., Boston,
MA.

\bibitem[Jouanolou(1980)]{Jou}
Jouanolou, J.~P., 1980.
\newblock{Id\'eaux r\'esultants.}
\newblock  Adv. in Math.  37, no. 3, 212--238.

\bibitem[Macaulay(1902)]{Mac}
Macaulay, F., 1902.
\newblock{Some formulae in elimination.}
\newblock Proc. London Math. Soc. 1 \textbf{33}, 3--27.


\bibitem[Macaulay(1916)]{Mac2}
Macaulay, F., 1916.
\newblock{The algebraic theory of modular systems.}
\newblock Cambridge University Press.

\bibitem[Pedersen \& Sturmfels(1996)]{PS}
Pedersen, P., Sturmfels, B., 1996.
\newblock{Mixed monomial bases.}
\newblock Algorithms in algebraic geometry and applications (Santander, 1994), 307--316,
Progr. Math., 143, Birkhäuser, Basel.

\bibitem[van der Waerden(1950)]{VdW}
van der Waerden, B.~L. 1950.
\newblock{Modern Algebra.}
\newblock $3^{rd}$ edn. New York, F.~Ungar Publishing Co.

\bibitem[Szanto(2002)]{sza}
Szanto,  A., 2002.
\newblock{Multivariate subresultants using Jouanolou's resultant matrices.}
\newblock Preprint.
\end{thebibliography}
\end{document}